\documentclass[a4paper, 12pt]{article}

\usepackage[koi8-r,cp1251]{inputenc}
\usepackage[english,russian]{babel}
\usepackage{amsfonts,amssymb,amsmath, hyperref}
\usepackage[final]{epsfig}
\usepackage{graphicx}

\begin{document}

\begin{flushleft}
 УДК 517.9
\end{flushleft}

\begin{center}
{\bf НЕСКОЛЬКО УТВЕРЖДЕНИЙ, РАВНОСИЛЬНЫХ ГИПОТЕЗЕ РИМАНА}
\end{center}

\begin{center}
  {\bf А. Р. Миротин}
  \end{center}

\begin{center}
\textit{Гомельский государственный университет им. Ф Скорины,
Гомель}
\end{center}

\begin{center}
{\bf  SOME ASSERTIONS WHICH ARE EQUIVALENT TO RIEMANN HYPOTHESIS}
\end{center}

\begin{center}
{\bf A.R. Mirotin}
\end{center}

\begin{center}
\textit{F. Scoryna Gomel State University, Gomel}
\end{center}

Сформулировано несколько утверждений, относящихся к гармоническому
анализу на бесконечномерном торе, и доказана их равносильность
гипотезе Римана.

\bigskip
{\bf \textit{Ключевые слова}}: \textit{гипотеза Римана,
бесконечномерный тор, функция Мебиуса, функция Эйлера, формула
обращения.}

\bigskip
Some assertions in harmonic analysis on the infinite dimensional
torus are stated and their equivalence to  Riemann hypothesis is
proved.

\bigskip
\textbf{\textit{Key words}}: \textit{Riemann hypothesis, infinite
dimensional torus, M\"obius function, Eiler function, inverse
formula.}

\bigskip
 \textbf{\textit{Введение}}

\textit{Дзета-функция Римана} определяется для комплексных
$s=\sigma+it$ при $\sigma>1$ следующим образом:
$$
\zeta(s)=\sum\limits_{n=1}^\infty\frac{1}{n^s}.
$$
Для вещественных $s>1$ она рассматривалась еще Эйлером. В своей
эпохальной работе [1], опубликованной в 1859 г., Б. Риман
аналитически продолжил ее в область $s\ne 1$ и показал, что все
так называемые нетривиальные нули дзета-функции (т. е. нули,
отличные от  -2, -4, -6, ...) лежат в "критической полосе"\quad
$0\leq\sigma\leq 1$. Знаменитая \textit{гипотеза Римана} (далее
RH) утверждает, что все нетривиальные нули лежат на "критической
прямой"\quad $\sigma=1/2$. Эти исследования были предприняты Б.
Риманом в связи с гипотезой Лежандра-Гаусса, согласно которой
количество $\pi(x)$ простых чисел, не превосходящих $x$,
удовлетворяет соотношению $\pi(x):{\rm Li}(x)\to 1\ (x\to\infty)$,
где
$$
{\rm Li}(x):=\int\limits_2^x\frac{dt}{\log t}
$$
---  интегральный логарифм.  Б.~Риман связал это
утверждение с распределением комплексных нулей дзета-функции и
более того, дал точную формулу для $\pi(x)$, содержащую
нетривиальные нули дзета-функции. В упрощенном виде его формула
может быть записана следующим образом (если число $x$ не есть
степень простого):
$$
\sum\limits_{n\leq
x}\Lambda(n)=x-\sum\limits_\rho\frac{x^\rho}{\rho}-\log(2\pi)-\frac{1}{2}\log(1-x^{-2}),
$$
где $\Lambda(n)=\log p$, если $n$  есть степень простого, и
$\Lambda(n)=0$ в противном случае (функция Мангольдта), а сумма
$\sum_\rho\frac{x^\rho}{\rho}$ распространяется на все
нетривиальные нули $\rho$ дзета-функции и понимается как
$\lim_{T\to\infty}\sum_{|\rho|<T}$. При этом гипотеза
Лежандра-Гаусса равносильна утверждению $ \sum\limits_{n\leq
x}\Lambda(n)\sim x\ (x\to\infty)$, и из приведенной выше формулы
Римана легко следует, что эта гипотеза верна, если все
нетривиальные нули лежат слева от прямой $\sigma=1$.

 Следует отметить, что исследованиям Римана предшествовали замечательные
 работы П. Л. Чебышева по распределению простых чисел, в которых
 фигурировала дзета-функция вещественного переменного. Вкладом
 Б. Римана, который  трудно переоценить,  был  именно выход к комплексную область,
 позволивший в конечном
  счете  установить закон распределения простых в
натуральном ряде. А именно, в 1898 году Ж.~Адамар и Ш.~де ля
Валле-Пуссен   доказали гипотезу Лежандра-Гаусса, установив
отсутствие нулей дзета-функции на
 прямой $\sigma=1$. Эту теорему теперь называют теоремой о простых
 числах (асимптотическим законом распределения простых
 чисел в натуральном ряде).

 Важность гипотезы Римана  для теории чисел объясняется прежде
 всего тем, что ее справедливость равносильна справедливости
 теоремы о простых числах с наилучшим возможным остатком. Дело в том, что,
  как указал сам Б. Риман, RH равносильна следующему утверждению: $\forall\varepsilon>0$

 $$
\pi(x)={\rm Li}(x)+O(x^{1/2+\varepsilon}).
 $$

 Удивительным образом эта гипотеза связана также с рядом утверждений из
 таких разных разделов науки, как теория конечных групп
 преобразований, теория вероятностей,
 функциональный анализ    и даже квантовая механика. Например, доказано, что RH
 равносильна каждому из следующих утверждений:

 --- $\sigma(n)<H_n+e^{H_n}\log H_n$,

 где  $\sigma(n)$ есть число делителей числа $n, H_n=1+1/2+\dots+1/n$ (Лагариас, 2002);

---  наибольший порядок $g(n)$ элементов симметрической группы $S_n$ при
достаточно больших $n$ удовлетворяет неравенству
$$
\log g(n)< \frac{1}{\sqrt{{\rm Li}(n)}};
$$

(Массиас, Николас, Робин,  1989);

--- для любой ненулевой комплексной функции
$g\in C^\infty_0(0,\infty)$ справедливо неравенство
$$
\sum\limits_\rho g^*(\rho)\overline{g^*(1-\rho)}>0,
$$
где
$$
g^*(s)=\int\limits_0^\infty g(x)x^{s-1}dx
$$
--- преобразование Меллина функции $g$ (А. Вейль, Э. Бомбьери);

--- для всех натуральных $n$
$$
\sum\limits_\rho(1-(1-1/\rho)^n)\geq 0
$$
(Х.-Дж. Ли, 1998);

--- при ${\rm Re}s>1/2$
$$
{\rm Re}\left(\frac{\xi'(s)}{\xi(s)}\right)>0,
$$
где
$\xi(s)=\frac{1}{2}s(s-1)\pi^{-\frac{s}{2}}\Gamma(\frac{s}{2})\zeta(s)$
(Хинкканен, 1998);

---
$$
\int\limits_0^\infty\int\limits_{1/2}^\infty\frac{1-12y^2}{(1+4y^2)^3}\log|\zeta(x+iy)|dxdy
=\pi\frac{3-\gamma}{32},
$$
где $\gamma$ --- постоянная Эйлера (Волчаков, 1998);

--- линейная оболочка семейства функций $f_a(t)=\{\frac{a}{t}\}-a\{\frac{1}{t}\},\ a\in (0;1)$
(фигурные скобки обозначают дробную часть числа) плотна в
$L^2[0;1]$ (Бёрлинг, Ниман).

Еще Б. Риман установил функциональное уравнение для дзета-функции,
которое можно записать в виде $\xi(s)=\xi(1-s)$. Но доказано, что
дзета-функция не удовлетворяет никакому дифференциальному
уравнению, что считается причиной столь сложного распределения ее
нулей.

 По-видимому, сейчас RH ---
 самая знаменитая нерешенная математическая проблема.
Следует отметить, что в настоящий момент мало кто сомневается в
справедливости этой гипотезы, поскольку многочисленные факты
свидетельствуют в ее пользу. Например, доказано, что первые 1,5
миллиарда комплексных нулей дзета-функции (расположенных в порядке
возрастания мнимых частей) лежат на критической прямой (ван де
Люн, Риэль, Винтнер, 1986). Там же лежат и $3\times 10^8$ ее нулей
с мнимыми частями из промежутка $[0;2\times 10^{20}]$, а также все
нули с мнимыми частями из промежутка $[10^{22};2\times
10^{22}+10^{10}]$ (А. Олдыжко). Известно также, что на критической
прямой расположено более 40 процентов комплексных нулей
дзета-функции (Левинсон, Сельберг, Конри). Доказано, что,
гипотетические исключения из RH могут быть расположены на
комплексной плоскости крайне редко (Бор, Ландау, Карлсон, Ингам).
Доказаны аналоги RH, относящиеся к дзета-функциям алгебраических
многообразий над конечными полями (А. Вейль, П. Делинь), некоторые
следствия RH получили независимые доказательства и т. д. Кроме
того, опровержение RH внесло бы хаос в распределение простых
чисел, и считается маловероятным, что природа может быть столь
извращенной.

 Далее нам понадобятся равносильные RH
 утверждения, относящиеся к  функциям Мебиуса $\mu$ и Эйлера $\varphi$.
Напомним, что по   определению $\mu(N_k)=(-1)^k$ где $N_k=p_1\dots
p_k$ есть произведение первых $k$ простых чисел, и $\mu(n)=0$ для
остальных натуральных $n$; $\varphi(n)$ есть число натуральных
чисел, меньших $n$ и взаимно простых с $n$.

\textbf{\textit{Теорема} 1.} (Дж. Литтлвуд) \textit{RH равносильна
каждому из следующих
 утверждений}:

 (i) $\forall \varepsilon>0\sum_{n\leq x}
\mu(n)=O(x^{1/2+\varepsilon})$;

 (ii) $\exists A>0\sum_{n\leq x}
\mu(n)=O(x^{1/2}\exp(A\log x/\log\log x))$.

\textbf{\textit{Теорема} 2}. (Николас, 1985) \textit{RH
равносильна каждому из следующих
 утверждений}:  $\forall k/\forall k,\mbox{ \textit{кроме конечного числа}},$

$$
\varphi(N_k)<\frac{e^{-\gamma}N_k}{\log\log N_k},
$$
\textit{где $N_k=p_1\dots p_k$ есть произведение первых $k$
простых чисел.}

Ниже с помощью гармонического анализа на полугруппах мы выведем из
этих результатов несколько новых утверждений, равносильных RH.
Всюду далее систематически будут использоваться следующие
обозначения: если $n=p_1^{\alpha_1(n)}p_2^{\alpha_2(n)}\dots$
--- каноническое разложение натурального числа $n$
 на простые множители, то $\alpha(n):=(\alpha_1(n), \alpha_2(n),\dots)$ ---
 финитная последовательность,
 и для $t\in {\Bbb T}^\omega$ мы полагаем
 $t^{\alpha(n)}:=t_1^{\alpha_1(n)}t_2^{\alpha_2(n)}\dots$
 (${\Bbb T}^\omega$
  обозначает бесконечномерный тор, т. е. счетное произведение единичных окружностей). Тогда, например,
  из нижеследующего предложения 2 вытекает, что RH равносильна
  такому утверждению:

$\forall k/\forall k,\mbox{ \textit{кроме конечного числа,
справедливо неравенство}}$
$$
\int\limits_{{\Bbb T}^\omega}\frac{g(t)}{t_1\dots
t_k}dt<\frac{1}{e^\gamma N_k^2\log\log N_k},\eqno(1)
$$
\textit{где} $g(t)=\sum_{n\geq
1}\frac{t^{\alpha(n)}}{n^2}/\sum_{n\geq
1}\frac{t^{\alpha(n)}}{n^3}$, а $dt$ --- нормированная мера Хаара
группы ${\Bbb T}^\omega$.

Более подробную информацию, касающуюся   гипотезы Римана, можно
найти в [2] --- [6].

\bigskip
\textbf{\textit{1 Вспомогательные сведения}}

В этом разделе будут изложены сведения по гармоническому анализу
на полугруппах, необходимые для дальнейшего. Систематическое
изложение этих (и других) вопросов гармонического анализа
содержится в [7] --- [10].

Пусть $S$ --- абелева полугруппа с сокращениями и нейтральным
элементом, записываемая мультипликативно, $G=S^{-1}S$ --- группа
ее частных.  Через $\widehat S$ обозначим мультипликативную
полугруппу всех \textit{ограниченных полухарактеров } полугруппы
$S$ (т.~е. ненулевых гомоморфизмов из $S$ в замкнутый единичный
диск $\overline{\Bbb{D}}$ комплексной плоскости с операцией
умножения), наделенную топологией поточечной сходимости
(превращающей ее в компактную топологическую полугруппу), а через
$\widehat S_+$ --- ее компактную подполугруппу, состоящую из
неотрицательных полухарактеров. Для $\rho\in\widehat S_+$ положим
также
$$
\widehat S_\rho:=\{\psi\in\widehat S:|\psi|\leq\rho\}.
$$

 \textit{Характером} полугруппы $S$ будем называть полухарактер,
равный по модулю единице; (компактная топологическая) группа
характеров полугруппы $S$ будет обозначаться $X$.

Необходимые нам сведения из [7] мы изложим в виде нескольких лемм.

\textbf{\textit{Лемма} 1} (полярное разложение полухарактера).
\textit{Любой полухарактер $\psi\in \widehat S$ можно представить
в виде
$$
\psi=\rho\chi,
$$
где $\rho\in\widehat S_+$, а $\chi\in X$.}

Для полухарактера $\rho\in\widehat S_+$ через $l_1(\rho)$ мы
обозначим нормированное пространство тех комплекснозначных функций
$f$ на группе $G$ с не более чем счетным носителем, которые
сосредоточены на $S$, и для которых
$$
\|f\|:=\sum\limits_{s\in S}|f(s)|\rho(s)<\infty.
$$

\textbf{\textit{Лемма} 2}. \textit{Пространство $l_1(\rho)$ есть
унитальная полупростая коммутативная банахова алгебра со сверткой
$$
f\ast g(s):=\sum\limits_{xy=s}f(x)g(y)
$$
в качестве умножения.}

Пусть $f\in l_1(\rho)$. Функцию $\tilde f$ на $\widehat S_\rho$,
определенную равенством
$$
\tilde f(\psi)=\sum\limits_{s\in S} f(s)\overline{\psi(s)} ,
$$
\noindent будем называть {\it преобразованием Лапласа}
 функции $f$.

\textbf{\textit{Лемма} 3}. \textit{Комплексные гомоморфизмы
алгебры $l_1(\rho)$ имеют в точности вид $f\mapsto \tilde f(\psi)
(\psi\in \widehat S_\rho)$. В частности, $\widetilde{f\ast
g}=\widetilde f\widetilde g$ при $f, g\in l_1(\rho)$.}

 Пусть $\psi=\rho\chi$ --- полярное разложение полухарактера $\psi$.
 Тогда ясно, что $\tilde f(\psi)=\widehat{(f\rho)}(\chi)$,
 где "крышка"\quad обозначает преобразование Фурье на группе $G$.
 Формула обращения для преобразования Фурье влечет теперь
 следующее утверждение.

\textbf{\textit{Лемма} 4} (формула обращения для преобразования
Лапласа). \textit{Пусть $f\in l_1(\rho)$. Если функция
$\chi\mapsto \widehat{(f\rho)}(\chi)$ принадлежит $L^1(X)$, то при
всех  $s$, для которых $\rho(s)\ne 0$, справедливо равенство
$$
f(s)=\frac{1}{\rho(s)}\int\limits_X\tilde f(\rho\chi)\chi(s)d\chi,
$$
где $d\chi$ --- нормированная мера Хаара группы} $X$.

\bigskip
\textbf{\textit{2 Формулировка и доказательство основных
результатов}}

Теперь мы в состоянии установить наши основные результаты.

\textbf{\textit{Предложение} 1}. \textit{RH равносильна каждому из
следующих утверждений}:

(i) \textit{для любого $\varepsilon>0$ и любого/некоторого
$\beta>1$ справедлива оценка}
$$
\int\limits_{\Bbb{T}^{\omega}}\sum\limits_{a\leq x}
\frac{1}{\sum\limits_{n\geq
1}(an)^{-\beta}t^{\alpha(an)}}dt=O(x^{1/2+\varepsilon});
$$

(ii) \textit{существует такое $A>0$ что для любого/некоторого
$\beta>1$ справедлива оценка}
$$
\int\limits_{\Bbb{T}^{\omega}}\sum\limits_{a\leq x}
\frac{1}{\sum\limits_{n\geq
1}(an)^{-\beta}t^{\alpha(an)}}dt=O(x^{1/2}\exp(A\log x/\log\log
x)).
$$

 \textit{Доказательство.}
Пусть $S$ есть мультипликативная полугруппа $\Bbb{N}^*$
натуральных чисел. Поскольку это свободная абелева полугруппа,
системой образующих которой служит множество простых чисел
$\{p_1,p_2,\ldots \}$, каждый ограниченный полухарактер полугруппы
$S$ имеет вид
$$
\psi(a)=\overline z^{\alpha(a)},
$$
 где
$z=(z_1,z_2,\ldots)\in\overline{\mathbb{D}}^{\infty}$,
$z^{\alpha(a)}=z_1^{\alpha_1(a)}z_2^{\alpha_2(a)}\ldots$,
$z_j=\overline{\psi(p_j)}$,
$a=p_1^{\alpha_1(a)}p_2^{\alpha_2(a)}\ldots$--- каноническое
разложение числа $a$ на простые множители. Таким образом, мы можем
считать, что $\widehat{\mathbb{N}^*}_+=[0,1]^{\omega}$,
$\widehat{\mathbb{N}^*}=\overline{\mathbb{D}}^{\omega}$, а группа
характеров полугруппы $\mathbb{N}^*$ есть  $\mathbb{T}^{\omega}$ (через
$K^\omega$ мы обозначаем тихоновское произведение счетного числа
экземпляров компакта $K$).

 Для этого случая, отождествляя полухарактер $\psi(a)=\overline
z^{\alpha(a)}$ с соответствующей точкой
$z\in\overline{\Bbb{D}}^{\omega}$, получаем, что преобразование
Лапласа имеет вид
$$
\tilde f(z)=\sum\limits_{a=1}^\infty f(a)z^{\alpha(a)}.
$$

Если полухарактеру $\rho\in \widehat{\Bbb{N}^*}_+$ соответствует
точка $r$ из $[0,1]^{\omega}$, то алгебра $l_1(\rho)$ состоит из
всех арифметических функций  $f$, для которых абсолютно сходится
ряд $\sum_{a=1}^\infty f(a)r^{\alpha(a)}$. Свертка в $l_1(\rho)$
задается формулой
$$
(f\ast g)(a)=\sum\limits_{d|a}f(d)g(\frac{a}{d}),
$$
\noindent  где суммирование распространяется на все положительные
делители числа $a$ ("свертка Дирихле"); функция $1_{\{1\}}$
(индикатор одноточечного множества $\{1\}$) будет единицей этой
алгебры.

В этом случае, при условии, что функция $t\mapsto \tilde f(r.t)$
принадлежит $L^1(\Bbb{T}^{\omega}, dt)$, формула обращения для
функций $f\in l_1(\rho)$ ($\rho(a)=r^{\alpha(a)},\ r\in
[0;1]^\omega$) принимает вид
 $$
f(a)=\frac{1}{r^{\alpha(a)}}\int\limits_{\Bbb{T}^{\omega}}\tilde
f(r.t)t^{-\alpha(a)}dt,\eqno(2)
 $$
\noindent где $dt$ --- нормированная мера Хаара группы
$\Bbb{T}^{\omega}$,  а $r.t$ обозначает последовательность
$(r_1t_1,r_2t_2,\ldots )$.

Выберем теперь $\rho(a)=a^{-\beta}\ (\beta>1$ фиксировано),  т. е.
положим $r=(p_1^{-\beta},p_2^{-\beta},\ldots)$. Очевидно, что
тогда $\mu\in l_1(\rho)$.  Применяя к известному тождеству
$$
\mu\ast 1=1_{\{1\}}
$$
преобразование Лапласа, получаем
  при $|z_j|<p_j^{-\beta},\ j\in\Bbb{N}$
$$
\tilde{\mu}(z)=\frac{1}{\sum_{a=1}^\infty z^{\alpha(a)}}.
$$

Покажем, что условия, достаточные для справедливости формулы
обращения, здесь выполнены при $\beta>1$, т. е. что функция
$t\mapsto \tilde \mu(r.t)$ принадлежит $L^1(\Bbb{T}^{\omega},
dt)$. В самом деле, эта функция есть преобразование Фурье функции,
определенной на мультипликативной группе $\mathbb{Q}_+$
положительных рациональных чисел (являющейся группой  частных
полугруппы $\Bbb{N}^*$), но сосредоточенной на полугруппе
$\mathbb{N}^*$ и равной на ней $(\mu\rho)(n)=\mu(n)n^{-\beta}$.
Поскольку эта функция, очевидно, принадлежит $L^2(\mathbb{Q}_+)$,
то функция $t\mapsto \tilde \mu(r.t)$ принадлежит
$L^2(\Bbb{T}^{\omega})\subset L^1(\Bbb{T}^{\omega})$ по теореме
Планшереля.

 Значит, для функции Мебиуса имеем при любом $\beta>1$ в
силу формулы обращения следующее интегральное представление (у нас
$r^{\alpha(n)}=n^{-\beta}$):
$$
\mu(a)=a^{\beta}\int\limits_{\Bbb{T}^{\omega}}\frac{1}{\sum_{n\geq
1}n^{-\beta}t^{\alpha(n)}}t^{-\alpha(a)}dt.\eqno(3)
$$

Следовательно,
$$
\sum\limits_{a\leq
x}\mu(a)=\int\limits_{\Bbb{T}^{\omega}}\sum\limits_{a\leq x}
\frac{1}{\sum\limits_{n\geq 1}(an)^{-\beta}t^{\alpha(an)}}dt
$$
и для завершения доказательства  осталось воспользоваться теоремой
1.

Более простой вид имеет следующий критерий справедливости RH.

\textbf{\textit{Предложение 2}}. \textit{RH равносильна каждому из
следующих утверждений}: $\forall k/\forall k,\mbox{ \textit{кроме
конечного числа, справедливо неравенство}}$
$$
\int\limits_{{\Bbb T}^\omega}\frac{g(t)}{t_1\dots
t_k}dt<\frac{1}{e^\gamma N_k^{\beta-1}\log\log N_k},
$$
\textit{где} $g(t)=\sum_{m\geq
1}\frac{t^{\alpha(m)}}{m^{\beta-1}}/\sum_{m\geq
1}\frac{t^{\alpha(m)}}{m^\beta},\ \beta>2$.

Доказательство. Будем рассуждать как в доказательстве предложения
1 (как и там, мы берем $\rho(a)=a^{-\beta}$, но считаем, что
$\beta>2$). Поскольку, очевидно, $\varphi(n)<n$, то $\varphi\in
l_1(\rho)$ . Применяя теперь преобразование Лапласа к известному
тождеству
$$
\varphi\ast 1={\rm id}
$$
 (здесь ${\rm id}(n)=n$),
получаем $\tilde{\varphi}(z)\cdot\tilde{1}(z)=\widetilde{\rm
id}(z)$ при $|z_j|<p_j^{-\beta},\ j\in\Bbb{N}$, откуда
$$
\tilde{\varphi}(z)=\frac{\sum_{m=1}^\infty
mz^{\alpha(m)}}{\sum_{m=1}^\infty z^{\alpha(m)}}.
$$

Покажем, что условия, достаточные для справедливости формулы
обращения, здесь выполнены при $\beta>2$. Действительно, функция
$t\mapsto \tilde \varphi(r.t)$ принадлежит $L^1(\Bbb{T}^{\omega},
dt)$, поскольку
$$
\left|\tilde \varphi(r.t)\right|=\frac{|\widetilde{\rm
id}(r.t)|}{|\tilde{1}(r.t)|}=|\widetilde{\rm id}(r.t)|\cdot|\tilde
\mu(r.t)|,
$$
причем первый сомножитель ограничен, $|\widetilde{\rm
id}(r.t)|\leq\sum_{m\geq 1}m^{1-\beta}<\infty$, а второй
принадлежит $L^1(\Bbb{T}^{\omega}, dt)$ (см. доказательство
предложения 1).

 Следовательно, по  формуле обращения (2) имеем при любом
$\beta>2$
 следующее интегральное представление функции Эйлера ($r^{\alpha(n)}=n^{-\beta}$):
$$
\varphi(n)=n^{\beta}\int\limits_{\Bbb{T}^{\omega}}\frac{\sum_{m=1}^\infty
m^{1-\beta}t^{\alpha(m)}}{\sum_{m=1}^\infty
m^{-\beta}t^{\alpha(m)}}t^{-\alpha(n)}dt, \eqno(4)
$$
и для завершения доказательства  осталось воспользоваться теоремой
2, так как  $t^{\alpha(N_k)}=t_1\dots t_k$.

Полагая $\beta=3$ в предложении 2, получаем утверждение (1),
сформулированное во введении.

 \textbf{\textit{Замечание 1}.} Теперь, когда формулы (3) и (4) получены
  с помощью гармонического анализа на полугруппах, можно дать их доказательство,
  являющееся концептуально более простым. Докажем формулу (3).
  Прежде всего заметим, что при $\beta>1$ справедливо следующее
  равенство:
  $$
\sum\limits_{n=1}^\infty\frac{t^{\alpha(n)}}{n^{\beta}}=
\sum\limits_{n=1}^\infty\frac{t_1^{\alpha_1(n)}t_2^{\alpha_2(n)}\dots}
{(p_1^\beta)^{\alpha_1(n)}(p_2^\beta)^{\alpha_2(n)}\dots}=
\prod\limits_{j=1}^\infty\sum\limits_{m=0}^\infty\left(\frac{t_j}{p_j^\beta}\right)^m=
\prod\limits_{j=1}^\infty\frac{1}{1-\frac{t_j}{p_j^\beta}}.\eqno(5)
$$

Поэтому для натурального $a$
$$
I(a):=\int\limits_{\Bbb{T}^{\omega}}\frac{1}{\sum\limits_{n=1}^\infty
n^{-\beta}t^{\alpha(n)}}t^{-\alpha(a)}dt=\int\limits_{\Bbb{T}^{\omega}}\frac{t^{-\alpha(a)}}
{\prod\limits_{j=1}^\infty\left(1-\frac{t_j}{p_j^\beta}\right)^{-1}}dt=
$$
$$
=\int\limits_{\Bbb{T}^{\omega}}\prod\limits_{j=1}^\infty\frac{1}
{\left(1-\frac{t_j}{p_j^\beta}\right)^{-1}t_j^{\alpha_j(a)}}dt=
\prod\limits_{j=1}^\infty\int\limits_{\Bbb{T}}\frac{1-\frac{t_j}{p_j^\beta}}{t_j^{\alpha_j(a)}}dt_j,
$$
где все меры $dt_j$ совпадают с нормированной мерой Хаара группы
$\Bbb{T}$ (последнее означает, что
$\int_{\Bbb{T}}f(t_j)dt_j=\int_0^{2\pi}f(e^{i\varphi})\frac{d\varphi}{2\pi}$
для любой непрерывной функции $f$ на $\Bbb{T}$).

Возможны два случая.

1) Число $a$ делится на квадрат простого, т. е. $\alpha_j(a)\geq
2$ при некотором $j$. В этом случае,  используя свойство
ортогональности характеров компактной группы (или прямым счетом),
получаем, что
$$
I_j(a):=\int\limits_{\Bbb{T}}\frac{1-\frac{t_j}{p_j^\beta}}{t_j^{\alpha_j(a)}}dt_j=0,
$$
а потому  $I(a)=0$, что доказывает (3) для выбранных $a$.

2)  Число $a$ свободно от квадратов, т. е. $\alpha_j(a)\in\{0;1\}$
при всех $j$. Как и в предыдущем случае, здесь легко подсчитать,
что

$$
I_j(a)=\begin{cases}
1,  \mbox { если } \alpha_j(a)=0\\ 
-1/p_j^\beta,  \mbox { если } \alpha_j(a)=1
\end{cases}
$$

Таким образом, если $a=p_1\dots p_k$, то
$I(a)=\prod_jI_j(a)=(-1)^k/a^\beta$, что равносильно формуле (3) в
этом случае.

Докажем формулу (4). Из формулы (5) следует, что для натурального
$a$
$$
J(a):=\int\limits_{\Bbb{T}^{\omega}}\frac{\sum_{n\geq 1}
n^{1-\beta}t^{\alpha(n)}}{\sum_{n\geq 1}
n^{-\beta}t^{\alpha(n)}}t^{-\alpha(a)}dt=
\prod\limits_{j=1}^\infty\int\limits_{\Bbb{T}}\frac{1-\frac{t_j}{p_j^\beta}}
{\left(1-\frac{t_j}{p_j^{\beta-1}}\right)t_j^{\alpha_j(a)}}dt_j.
$$
Рассмотрим интеграл
$$
J_j(a):=\int\limits_{\Bbb{T}}\frac{1-\frac{t_j}
{p_j^\beta}}{\left(1-\frac{t_j}{p_j^{\beta-1}}\right)t_j^{\alpha_j(a)}}dt_j=
\frac{1}{p_j}\int\limits_{\Bbb{T}}\frac{p_j^\beta-t_j}
{(p_j^{\beta-1}-t_j)t_j^{\alpha_j(a)}}dt_j.
$$
Поскольку
$$
\frac{p_j^\beta-t_j}
{p_j^{\beta-1}-t_j}=1+(p_j-1)\sum\limits_{m=0}^\infty\left(\frac{t_j}{p_j^{\beta-1}}\right)^m,
$$
то, пользуясь ортонормированностью характеров группы $\Bbb{T}$,
получим
$$
J_j(a)=
\begin{cases}1,  \mbox { если } \alpha_j(a)=0 \\
\frac{p_j-1}{p_j^{(\beta-1)\alpha_j(a)+1}},  \mbox { если }
\alpha_j(a)\geq 1
\end{cases}.
$$
Поэтому, если $a=p_1^{\alpha_1(a)}\dots p_k^{\alpha_k(a)}$, то
$$
J(a)=\prod_jJ_j(a)=\prod\limits_{j=1}^k\frac{p_j-1}{p_j^{(\beta-1)\alpha_j(a)+1}}
=\frac{\varphi(a)}{a^\beta}
$$
 (мы воспользовались тем, что
$\varphi(a)=a\prod_{j=1}^k(1-1/p_j)$). Последнее равенство
равносильно формуле (4).

\textbf{\textit{Замечание 2.}} Подобно формулам (3) и (4) могут
быть получены интегральные представления и других арифметических
функций. Так, например, исходя из  тождества
$$
\Lambda\ast 1=\log
$$
(см. [11], с. 145), для функции Мангольдта  $\Lambda$ имеем при
$\beta>1$
$$
\Lambda(n)=n^{\beta}\int\limits_{\Bbb{T}^{\omega}}\frac{\sum_{m=1}^\infty
 m^{-\beta}\log m t^{\alpha(m)}}{\sum_{m=1}^\infty
m^{-\beta}t^{\alpha(m)}}t^{-\alpha(n)}dt.
$$

\bigskip
\begin{center}
  ЛИТЕРАТУРА
  \end{center}

1. \textit{Riemann, B.}  Ueber die Anzahl der Primzahlen unter
einer gegebenen Gr\"{o}sse / B. Riemann - in: Monat. der
K\"{o}nigl.  ---  Berlin : Preuss. Akad. der Wissen., 1859 (1860).
- P. 671 -- 680.

2.\textit{ Corney, J. B.} The Riemann Hypothesis / J. B. Corney //
Notices of the AMS. -- 2004. -- Vol. 50, N 3. - P. 341 -- 353.

3. \textit{Bombieri, E.} The Riemann Hypothesis / E. Bombieri ---
in.: The Millennium Prize Problems. --- Providence, RI : AMS,
2006. -- P. 107 --  129.

4. The Riemann Hypothesis / http://www.aimath.org/WWN/rh/
[Электронный ресурс].

5. \textit{Titchmarsh, E. S.} The Theory of the Riemann Zeta
Function / E. S. Titchmarsh --- 2nd ed. revised by R. D.
Heath-Brown, ---  Oxford :  Oxford  University Press, 1986.

6.\textit{ Ivi\v{c}, A.} The  Riemann Zeta-Function --- The Theory
of the Riemann Zeta-Function with Applications / A. Ivi\v{c} ---
New York : John Wiley, 1985.

7. \textit{Миротин, А. Р.} Гармонический анализ на абелевых
полугруппах / А. Р. Миротин. --- Гомель : ГГУ им. Ф. Скорины,
2008. -- 207 с.

8.\textit{ Rudin, W.} Fourier Analysis on Groups / W. Rudin ---
New York : Interscience Publishers, 1962. -- 285 p.

9.\textit{ Люмис, Л.} Введение в абстрактный гармонический анализ
/ Л. Люмис. --- М. : ИЛ, 1956. -- 251 с.

10. Mirotin A. R. Every Invariant Measure Semigroup Contains an Ideal
which is Embeddable in a Group. Semigroup Forum.
 1999. Vol.~59, No~3.  P. 354 -- 361.

11. \textit{Постников, А. Г. } Введение в аналитическую теорию
чисел / А. Г. Постников. --- М.: Наука, 1971. -- 416 с.

\end{document}